\documentclass[12pt,a4paper]{amsart}
\usepackage{amsmath,amssymb,amscd}
\input xy
\xyoption{all}
\newcommand{\ra}{\rightarrow}

\newcommand{\CC}{\mathbb C}
\newcommand{\ZZ}{\mathbb Z}

\newcommand{\HH}{\mathbb H}
\newcommand{\PP}{\mathbb P}
\newcommand{\cP}{\mathcal{P}}
\newcommand{\cE}{\mathcal{E}}
\newcommand{\cF}{\mathcal{F}}
\newcommand{\cU}{\mathcal{U}}
\newcommand{\cO}{\mathcal{O}}
\newcommand{\cM}{\mathcal{M}}

\theoremstyle{plain}
\newtheorem{theorem}{Theorem}[section]
\newtheorem{lemma}[theorem]{Lemma}
\newtheorem{proposition}[theorem]{Proposition}
\newtheorem{cor}[theorem]{Corollary}
\newtheorem{remark}[theorem]{Remark}
\newtheorem{rem}[theorem]{Remark}
\newtheorem{ex}[theorem]{Example}

\begin{document}
\title[Coherent systems]{Coherent systems on elliptic curves}

\author{H. Lange}
\author{P. E. Newstead}

\address{H. Lange\\Mathematisches Institut\\
              Universit\"at Erlangen-N\"urnberg\\
              Bismarckstra\ss e $1\frac{ 1}{2}$\\
              D-$91054$ Erlangen\\
              Germany}
              \email{lange@mi.uni-erlangen.de}
\address{P.E. Newstead\\Department of Mathematical Sciences\\
              University of Liverpool\\
              Peach Street, Liverpool L69 7ZL, UK}
\email{newstead@liv.ac.uk}
\thanks{Both authors are members of the research group VBAC (Vector Bundles on Algebraic Curves). The second author 
         acknowledges support from EPSRC Grant No. EP/C515064, and would like to thank the Mathematisches Institut der Universit\"at 
         Erlangen-N\"urnberg for its hospitality}

\begin{abstract}
In this paper we consider coherent systems $(E,V)$ on an elliptic curve which 
are ${\alpha}$-stable with respect to some value of a parameter 
$\alpha$. We show that the corresponding moduli spaces, if non-empty, are smooth and irreducible of the expected dimenson.
Moreover we give precise conditions for non-emptiness of the moduli spaces. Finally we study the variation of the moduli spaces with $\alpha$.
\end{abstract}
\maketitle

\section{Introduction}

\noindent
A {\it coherent system of type $(n,d,k)$} on a smooth projective curve $X$ over an algebraically closed field is by 
definition a pair $(E,V)$ consisting of
a vector bundle $E$ of rank $n$ and degree $d$ over $X$ and a vector subspace $V \subset H^0(E)$ of dimension $k$. 
For any real number $\alpha$, the {\it $\alpha$-slope} of a coherent system $(E,V)$ of type $(n,d,k)$ is defined by
$$
\mu_{\alpha}(E,V) := \frac{d}{n} + \alpha \frac{k}{n}.
$$
A {\it coherent subsystem} of $(E,V)$ is a coherent system $(E',V')$ such that $E'$ is a subbundle of $E$ and 
$V' \subset V \cap H^0(E')$.
A coherent system $(E,V)$ is called 
{\it $\alpha$-stable} ({\it $\alpha$-semistable}) if
$$
\mu_{\alpha}(E',V') < \mu_{\alpha}(E,V) \ \ (\mu_{\alpha}(E',V') \le \mu_{\alpha}(E,V))
$$
for every proper coherent subsystem $(E',V')$ of $(E,V)$.
The $\alpha$-stable coherent systems of type $(n,d,k)$ on $X$ form a quasiprojective moduli space which we 
denote by $G(\alpha;n,d,k)$.

The general theory of \cite{bgn} is true for coherent systems on curves of all genera, but most of the detailed results require 
$g \geq 2$. In a previous paper \cite{ln} we studied the case $g=0$ and proved that all non-empty moduli spaces 
$G(\alpha;n,d,k)$ are smooth and irreducible of the expected dimension. Moreover we obtained results on non-emptiness for 
certain values of $k$ and $\alpha$.

In the current paper we study the case $g=1$. Our main result is that 
the non-empty moduli spaces $G(\alpha;n,d,k)$
are again smooth and irreducible of the expected dimension. 
We also determine precisely the conditions for non-emptiness.
Our results are therefore rather stronger than in the case $g=0$. 
The reasons for this are that there exist semistable bundles for any rank
$n$ and degree $d$ and that the general such bundle has the form 
$E=E_1 \oplus \cdots \oplus E_h$, where $h=\gcd(n,d)$ and the $E_i$ 
are stable of the same slope; these facts are 
due to Atiyah \cite{at} and Tu \cite{tu}.

In order to state our main results more precisely, we first recall the
definition of the {\em Brill-Noether number} $\beta(n,d,k)$ 
\cite[Definition 2.7]{bgn}. For an elliptic curve, this is independent of
$n$ and we write
$$
\beta(d,k):=\beta(n,d,k)=k(d-k)+1.
$$
We write also $M(n,d)$ for the moduli space of stable bundles of rank $n$ 
and degree $d$ on $X$. Note that, if $k\ge1$ and $\alpha\le0$, there do not 
exist $\alpha$-stable coherent systems of type $(n,d,k)$ (see 
\cite[section 2.1]{bgn}). Our main result can now be summarised as follows.

\noindent{\bf Theorem.}\begin{em} Let $X$ be an elliptic curve and suppose
$n\ge1$, $k\ge0$. Then
\begin{itemize}
\item[(i)] if $G(\alpha;n,d,k)\ne\emptyset$, it is smooth and irreducible 
of dimension $\beta(d,k)$;
\item[(ii)] $G(\alpha;n,d,0)\simeq M(n,d)$ for all $\alpha$; in particular
it is non-empty if and only if $\gcd(n,d)=1$;
\item[(iii)] for $\alpha>0$ and $k\ge1$, $G(\alpha;1,d,k)$ is independent of
$\alpha$ and is non-empty if and only if either $d=0$, $k=1$ or $k\le d$;
\item[(iv)] for $\alpha>0$, $n\ge2$ and $k\ge1$, $G(\alpha;n,d,k)\ne\emptyset$
if and only if $(n-k)\alpha<d$ and either $k<d$ or $k=d$ and $\gcd(n,d)=1$.
\end{itemize}\end{em} 

Following some preliminaries in section 2, we prove (ii) and (iii) (together
with (i) for the cases $k=0$ and $n=1$) in section 3. In section 4, we prove
(i) (Theorem \ref{thm3.3}) and  the necessity of the condition in (iv)
(Corollary \ref{cor3.2} and Proposition \ref{prop3.4}). In this section, we show also
that the birational type of $G(\alpha;n,d,k)$ is independent of $\alpha$
(Theorem \ref{cor3.4}) and give descriptions of the generic elements of
$G(\alpha;n,d,k)$ (Propositions \ref{prop3.5}, \ref{prop3.6} and \ref{prop3.7}).
In section 5, we prove the sufficiency of the condition in (iv) (Theorems
\ref{thm4.1}, \ref{thm4.2}
and \ref{thm4.3}). In section 6,
we consider the variation of the spaces $G(\alpha;n,d,k)$ as $\alpha$ varies. We discuss the 
``flips'' introduced in \cite{bgn} and give examples of cases in which 
there are no flips and some in which flips genuinely exist. 

In order to show that the small-$\alpha$ and large-$\alpha$ moduli spaces can be non-isomorphic varieties, we discuss in 
section 7 the Poincar\'e and Picard bundles for elliptic curves. These bundles have been much studied for general genus, 
but we could not find an account in the literature for the elliptic case. 
In section 8 we use a calculation of the 
Chern classes to give an example in which small-$\alpha$ and large-$\alpha$ moduli spaces are non-isomorphic although they are both Grassmannian bundles over
$X$ of the same dimension.

We refer to \cite{bgn} for general information about coherent systems on
algebraic curves and additional references (see also \cite{he}, where much of
the general theory of coherent systems is developed).

\noindent{\bf Acknowledgment.} We would like to thank Montserrat Teixidor
i Bigas for drawing our attention to the reference \cite{bdw}, and the 
referee for a thorough reading of the paper and some useful comments which
have led to improvements in the presentation.

\section{Preliminaries}

\noindent
Let $X$ be an elliptic curve which for simplicity we suppose to be defined over an algebraically closed field of characteristic 0.
We begin by recalling some facts on vector bundles on $X$.

A complete classification of indecomposable bundles on $X$ is given in 
\cite{at}; in particular, the indecomposable bundles of any fixed rank and
degree form a family parametrised by $X$ \cite[Theorem 7]{at}. 
Since every vector bundle can be written as 
a direct sum of indecomposable ones, uniquely up to order of direct summands, 
this gives a complete description of all vector bundles on $X$. 

We have the following facts covering the relationship between
indecomposability and stability:
\begin{itemize}
\item every indecomposable bundle on $X$ is semistable,
\item an indecomposable bundle is stable if and only its 
rank and degree are coprime, 
\item every simple bundle is stable. 
\end{itemize}
These are consequences of \cite{at}. For detailed proofs of 
the first two facts, see \cite[Appendix A]{tu}; for the third, we 
need only show that, if $E$ is simple, then its rank and degree are coprime.
This last follows from \cite[Theorem 5 and Lemmas 24, 26]{at}.

It follows from the second fact above that, if $\gcd(n,d)=1$, then the 
moduli space $M(n,d)$ is isomorphic to $X$.
In \cite{tu} it is shown more generally
that the moduli space of S-equivalence classes of semistable bundles of rank $n$ and degree $d$ on $X$ is isomorphic to the 
$h$-th symmetric product $S^hX$ of $X$ where $h = \mbox{gcd}(n,d)$. 
In fact, every point of $S^hX$ is represented by a unique polystable bundle
$$
E = E_1 \oplus \cdots \oplus E_h,
$$
where each $E_i$ is stable of rank $\frac{n}{h}$ and degree $\frac{d}{h}$.

We need from \cite{at} the description of the indecomposable vector bundles of degree 0 on $X$. Every such bundle can be 
written uniquely in the form $F_r \otimes L$ with $L \in \mbox{Pic}^0(X)$, where $F_r$ is defined inductively 
as follows:
\begin{itemize}
\item $F_1=\cO$, 
\item $F_r$ is determined uniquely as a nontrivial extension
\begin{equation} \label{eqn1.1}
0 \ra \cO \ra F_r \ra F_{r-1} \ra 0.
\end{equation}
\end{itemize}
Moreover 
\begin{equation} \label{eqn1.2}
h^0(F_r) = h^1(F_r) = 1
\end{equation}
and, for any nontrivial $L \in \mbox{Pic}^0(X)$, 
$$
h^0(F_r \otimes L) = h^1(F_r \otimes L) = 0.
$$
For an indecomposable bundle $E$ of positive degree,
\begin{equation} \label{eqn1.3}
h^1(E) = 0 \quad \mbox{and} \quad h^0(E) = \deg (E).
\end{equation}

We need the following two lemmas which are implicit in \cite{at}, but for which we could not find a direct reference.

\begin{lemma} \label{lemma2.1}
Let $E$ be an indecomposable vector bundle on $X$ of rank $n$ 
and degree $d$. If $d>n$, then $E$ is generated by 
its global sections.
\end{lemma}

\begin{proof}
Consider the exact sequence $0 \ra E(-p) \ra E \ra E_p \ra 0$ for any point $p \in X$.
Since $E(-p)$ is indecomposable of positive degree, we have $h^1(E(-p)) = 0$. This implies the assertion.
\end{proof}

\begin{lemma} \label{lemma2.2}
Suppose that $E = E_1 \oplus \cdots \oplus E_l$ with all $E_i$ indecomposable. Then $\dim$ {\em Aut}$(E) \geq l$. 
Moreover equality holds if and only if $E$ is polystable and the $E_i$ are pairwise non-isomorphic. 
\end{lemma}

\begin{proof}
The first statement being obvious, it remains to prove the second one. If $E$ is polystable and the $E_i$ are 
pairwise non-isomorphic, then each $E_i$ is simple and $\mbox{Hom}(E_i,E_j) = 0$ for $i \neq j$. So 
$$
\mbox{Aut}(E) =
\prod_{i=1}^l \mbox{Aut}(E_i) \simeq (\CC^*)^l.
$$
Conversely, suppose $\dim \mbox{Aut}(E) = l$. Then each $E_i$ is simple and thus stable. Moreover $\mbox{Hom}(E_i,E_j) = 0$
 for $i \neq j$. So the $E_i$ are pairwise non-isomorphic. If $\mu(E_i) < \mu(E_j)$, then $\deg(Hom(E_i,E_j)) > 0$ 
 and hence Hom$(E_i,E_j) \neq 0$, a contradiction. So all $E_i$ have the same slope and $E$ is polystable.  
\end{proof}

\section{The moduli spaces $G(\alpha;n,d,0)$ and $G(\alpha;1,d,k)$}

\noindent
Suppose first that $k=0$.

\begin{proposition}\label{propa1} $G(\alpha;n,d,0)\simeq M(n,d)$ for all
$\alpha$. In particular $G(\alpha;n,d,0)\ne\emptyset$ if and only if
$\gcd(n,d)=1$ and it is then smooth and irreducible of dimension
$\beta(d,0)=1$.
\end{proposition}
\begin{proof} It is immediate from the definitions that 
$G(\alpha;n,d,0)\simeq M(n,d)$. The rest of the proposition follows from the
results of \cite{at} and \cite{tu} (see section 2).
\end{proof}

We now suppose that $n=1$ and $k\ge1$.

\begin{proposition}\label{propa2}
For $\alpha>0$ and $k\ge1$, $G(\alpha;1,d,k)$ is independent of $\alpha$ and
is non-empty if and only if either $d=0$, $k=1$ or $k\le d$. It is then
smooth and irreducible of dimension $\beta(d,k)=k(d-k)+1$.
\end{proposition}
\begin{proof} When $n=1$, all coherent systems $(E,V)$ are $\alpha$-stable
for all $\alpha>0$. The existence result follows at once from the facts that 
${\mathcal O}$ 
is the unique line bundle of degree $0$ which has a non-zero section and that,
if $d>0$, then $h^0(L)=d$ for every line bundle $L$ of degree $d$ by
Riemann-Roch. It follows also that $G(\alpha;1,d,k)$ is a point if $d=0$,
$k=1$ and a Grassmannian bundle over $M(1,d)\simeq X$ with fibre $Gr(k,d)$
if $d>0$. This proves the remaining parts of the proposition.
\end{proof}

\begin{remark}\label{rmka} 
\begin{em} Of course, for $\alpha>0$, $G(\alpha;1,d,k)$ coincides
with the classical variety $G^{k-1}_d(X)$ of linear systems of degree $d$ 
and projective dimension $k-1$ on $X$ (see \cite[p.182]{acgh}). 
Proposition \ref{propa2} is therefore a classical result \cite[Chapter V]{acgh}.
\end{em}\end{remark}

\section{The moduli spaces $G(\alpha;n,d,k)$ for $n\ge2$}

\noindent
Since the canonical line bundle on $X$ is trivial, the {\it Petri map} of $(E,V)$ 
(given by multiplication of sections) is of the form
$$
V \otimes H^0(E^*) \ra H^0(E \otimes E^*).
$$
For any component $U$ of $G(\alpha;n,d,k)$ we have 
(see \cite[Corollaire 3.14]{he}) 
$$
\dim(U) \geq \beta(d,k).
$$
Moreover $G(\alpha;n,d,k)$ is smooth of the expected dimension 
$\beta(d,k)$ at $(E,V)$ if and only if the Petri map
is injective (see \cite[Proposition 3.10]{bgn}).

\begin{lemma} \label{lemma3.1}
Suppose $n\ge2$, $k>0$ and $(E,V)$ is an $\alpha$-stable coherent system 
on an elliptic curve $X$ of type $(n,d,k)$ for some $\alpha$. 
Then every indecomposable direct summand of $E$ is of positive degree.  
\end{lemma}

\begin{proof}
Let $F$ be an indecomposable direct summand of $E$ of degree $f\le0$. 
If $f<0$, we would have 
$(F,0)$ as a direct summand of the coherent system $(E,V)$, contradicting $\alpha$-stability for all $\alpha$.
If $f=0$ and the induced map $V \ra H^0(F)$ is 0, the same argument applies.

Otherwise $F= F_r$ for some $r$
(see section 2) and the map $V \ra H^0(F_r)$ is surjective.
In this case we can write
$$
E = F_r \oplus G.
$$
It follows from (\ref{eqn1.1}) and (\ref{eqn1.2}) that $F_r$ has a subbundle $\cO$ and that $H^0(\cO)$ maps isomorphically
to $H^0(F_r)$.
Since we are assuming $V \ra H^0(F_r)$ is nonzero, the 
coherent subsystem $(\cO,H^0(\cO))$ of $(F_r,H^0(F_r))$ lifts to a coherent subsystem of $(E,V)$. 
But $(E,V)$ also has a coherent subsystem $(G,V')$ of type $(n-r,d,k-1)$. The existence of these two coherent 
subsystems contradicts the $\alpha$-stability of $(E,V)$. 
\end{proof}

\begin{cor} \label{cor3.2}
Suppose $n\ge2$ and $k>0$. If $G(\alpha;n,d,k)$ is non-empty, then $k \le d$.
\end{cor}

\begin{proof}
For $(E,V) \in G(\alpha;n,d,k)$, the lemma implies that every indecomposable direct summand of $E$ is of positive degree.
Now, by (\ref{eqn1.3}), 
$$k=\dim V\le h^0(E)=d.
$$ \end{proof}

\begin{theorem} \label{thm3.3}
Suppose $n\ge2$, $k>0$ and $G(\alpha;n,d,k)$ is non-empty. 
Then $G(\alpha;n,d,k)$ is smooth and irreducible of dimension $\beta(d,k)$.
Moreover, for a general $(E,V) \in G(\alpha;n,d,k)$, $E$ is polystable. To be more precise,
$$
E = E_1 \oplus \cdots \oplus E_h
$$
with $h = \gcd(n,d)$ and all $E_i$ stable and pairwise non-isomorphic of the same slope.
\end{theorem}

\begin{proof}
If $(E,V)$ is $\alpha$-stable, if follows from Lemma \ref{lemma3.1} that $H^0(E^*) = 0$. Hence the Petri map 
$V \otimes H^0(E^*) \ra H^0(E \otimes E^*)$ is injective. So 
$G(\alpha;n,d,k)$ is smooth and has dimension $\beta(d,k)$.
By Lemma \ref{lemma3.1}, for every $(E,V) \in G(\alpha;n,d,k)$, $E$ is of the form
\begin{eqnarray}
E = E_1 \oplus \ldots \oplus E_l   \label{eqn3.1}
\end{eqnarray}
with $E_i$ indecomposable and $\deg (E_i) \geq 1$ for all $i$. 

Suppose $(n_1,d_1), \ldots , (n_l,d_l)$ are ordered pairs of positive integers with $\sum n_i = n$ and $\sum d_i = d$. Then 
the sequences $(E_1, \ldots, E_l)$ of indecomposable bundles with 
rk$(E_i) = n_i$ and $\deg (E_i) = d_i$ form a family parametrised 
by $X^l$ (see section 2). Note that $h^0(E_i) = d_i$. So $h^0(E_1 \oplus \cdots \oplus E_l) = d$. Therefore the set of all pairs 
$(E_1 \oplus \cdots \oplus E_l,V)$, where $\dim V = k$, forms a family 
parametrised by a projective variety $W$ (which is a
 Grassmannian bundle over $X^l$). Let $U$ denote the open subset of $W$ corresponding to $\alpha$-stable coherent systems.
 If $U \neq \emptyset$, there is a canonical morphism
 $$
 \varphi: U \ra G(\alpha;n,d,k).
 $$
 For any $(E,V) \in$ Im$\varphi$ we have
 $$
 \dim \varphi^{-1}(E,V) = \dim \mbox{Aut}(E) -1.
 $$
 So
$$
\begin{array}{ll}
\dim \mbox{Im}\varphi &= \dim Gr(k,h^0(E)) + l - \min \dim \mbox{Aut}(E) +1\\
         &= k(d-k) + l - \min \dim \mbox{Aut(E)} + 1\\
         &= \beta(d,k) +l - \min \dim \mbox{Aut}(E),
\end{array}
$$
where the minimum is to be taken over all $(E,V)$ with $E=E_1 \oplus \cdots \oplus E_l$ as above.
Hence, if $\min \dim \mbox{Aut}(E) > l$, this cannot give an open set of a component of $G(\alpha;n,d,k)$. 
By Lemma \ref{lemma2.2} this implies that, if the closure of Im$\varphi$ is an irreducible component of $G(\alpha;n,d,k)$,
then all $E_i$ are stable with the same slope. Hence $l = h= \mbox{gcd}(n,d)$ and $n_i = \frac{n}{h}, \;\;
d_i = \frac{d}{h}$. So $U$ is uniquely determined and $\overline{\mbox{Im}\varphi} = G(\alpha;n,d,k)$.
Moreover, for the general element of $U$, the bundles $E_i$ are pairwise non-isomorphic.              
\end{proof}

\begin{theorem} \label{cor3.4}
The set 
$$
I(n,d,k) := \{\alpha\;|\; G(\alpha;n,d,k) \neq \emptyset \}
$$
 is an open interval (possibly infinite or empty). Moreover
\begin{itemize}
\item[(i)] if $I(n,d,k)\ne\emptyset$, there exists a coherent system of type
$(n,d,k)$ which is $\alpha$-stable for all $\alpha\in I(n,d,k)$;
\item[(ii)] the birational type of $G(\alpha;n,d,k)$ is independent of
$\alpha\in I(n,d,k)$.
\end{itemize}\end{theorem}

\begin{proof}
Suppose $\alpha_1, \alpha_2 \in I(n,d,k)$ with $\alpha_1 < \alpha_2$. Then the irreducibility of the variety $W$ 
in the proof of Theorem \ref{thm3.3} shows that there exists a coherent system $(E_1 \oplus \cdots \oplus E_h,V) \in W$ 
which is simultaneously $\alpha_1$-stable and $\alpha_2$-stable. So $(E_1 \oplus \cdots \oplus E_h,V)$ is $\alpha$-stable for 
all $\alpha \in [\alpha_1,\alpha_2]$ by \cite[Lemma 3.14]{bg}. This proves that $I(n,d,k)$ is an interval. If $G(\alpha;n,d,k)
\neq \emptyset$ at an endpoint $\alpha_1$ of $I(n,d,k)$, it would be possible to extend the interval, a contradiction. This 
proves the openness of $I(n,d,k)$.

For (i), note that the general element of $W$ defines a coherent system
with the required property; (ii) follows at ance from (i) and the 
irreducibility of $G(\alpha;n,d,k)$.
\end{proof}

\begin{proposition} \label{prop3.4}
If {\em gcd}$(n,d) > 1$, then $G(\alpha;n,d,d)$ is empty for all $\alpha$.
\end{proposition}

\begin{proof} If $d=0$, this follows from Proposition \ref{propa1}. If $d>0$ and $G(\alpha;n,d,d) \neq \emptyset$, 
then, by Theorem \ref{thm3.3}, its general element is of the form 
$(E,V)$, where 
$E = E_1 \oplus \cdots \oplus E_h$ with $h \geq 2$ and $V = H^0(E)$. So $(E,V)$ is a direct sum of $h$ coherent systems and 
thus not $\alpha$-stable for any $\alpha$.
\end{proof}

\begin{proposition} \label{prop3.5}
Let $0 < k < n$ and suppose $G(\alpha;n,d,k)$ is non-empty. Then, for a general $(E,V) \in G(\alpha;n,d,k),$ 
we have an exact sequence
\begin{eqnarray}
0 \ra \cO^k \ra E \ra G \ra 0,    \label{eqn3.2} 
\end{eqnarray}
where $V = H^0(\cO^k) \subset H^0(E)$ and $G$ is a polystable vector bundle with pairwise non-isomorphic indecomposable direct 
summands. 
\end{proposition}

\begin{proof}
By Theorem \ref{thm3.3}, we may assume that $E$ is polystable with pairwise non-isomorphic indecomposable
direct summands. The subspace $V \subset H^0(E)$ is given by a homomorphism
$\cO^k \ra E$. We first claim that this homomorphism gives $\cO^k$ the structure of a subbundle of $E$, i.e. the homomorphism 
is injective and the quotient $G$ is a vector bundle.

Suppose this is not the case. Then there exists a global section of $E$, contained in $V$, with a zero. 
If $d \leq n$, this contradicts $\alpha$-stability.
If $d > n$, each indecomposable direct summand $E_i$ of $E$ is generated by its global sections according to Lemma \ref{lemma2.1}.
Hence, by \cite[Theorem 2]{at}, a general subspace $V$ embeds $\cO^k$ as a subbundle of $E$. This completes the proof of the claim.
   
It remains to show that $G$ is polystable. To see this, suppose
$$
G = G_1 \oplus \cdots \oplus G_l
$$
with indecomposable direct summands $G_i$. We have $\deg (G_i) >0$ for all $i$, since any quotient of a polystable 
bundle of positive degree is of positive degree. This implies $H^0(G^*) = 0$.
Let $\cM_l$ denote the moduli space of all sequences $(G_1, \ldots, G_l)$ for any fixed values of the ranks and degrees of the $G_i$.
Consider the subscheme
$U$ of $G(\alpha;n,d,k)$ parametrising the coherent systems given by exact sequences
$$
0 \ra \cO^k \ra E' \ra G_1 \oplus \cdots \oplus G_l \ra 0
$$
with $(G_1 \oplus \cdots \oplus G_l) \in \cM_l$ and $(E',H^0(\cO^k))$ $\alpha$-stable. Then either $U$ is empty or we have
$$
\begin{array}{ll}
&\dim U =\\
&= \dim \cM_l + \dim \mbox{Ext}^1(G,\cO^k) - \dim \mbox{Aut} (\cO^k) - \min \dim \mbox{Aut}(G) +1   \\
       & = l + kd - k^2 -\min \dim \mbox{Aut}(G) + 1\\
       & = \beta(d,k) + l - \min \dim \mbox{Aut}(G).
\end{array}
$$  
So, by Lemma \ref{lemma2.2}, $\dim U \leq \beta(d,k)$ with equality only if the $G_i$ are stable of the same slope.
The result follows as in the proof of Theorem \ref{thm3.3}.
\end{proof}

\begin{proposition} \label{prop3.6}
Suppose $n\ge2$ and $G(\alpha;n,d,n)$ is non-empty. For a general $(E,V) \in G(\alpha;n,d,n)$, there is an exact sequence
$$
0 \ra \cO^n \ra E \ra T \ra 0
$$
with $V = H^0(\cO^n) \subset H^0(E)$ and $T$ a torsion sheaf of length $d$.
\end{proposition}

\begin{proof}
By Corollary \ref{cor3.2} and Proposition \ref{prop3.4}, we have $d > n$. It follows from Lemma \ref{lemma2.1} and Theorem \ref{thm3.3}
that $E$ is generated by its sections. Hence, the general subspace of dimension $n$ of $H^0(E)$ generates $E$ generically.
This implies the assertion.
\end{proof}  

\begin{proposition} \label{prop3.7}
Suppose $k>n$ and $G(\alpha;n,d,k)$ is non-empty. For a general $(E,V) \in G(\alpha;n,d,k)$, there is an exact sequence
\begin{eqnarray}
0 \ra H \ra \cO^k \ra E \ra 0   \label{eqn3.3}
\end{eqnarray}
such that
\begin{itemize}
\item[(i)] $H^0(\cO^k) \ra H^0(E)$ is an isomorphism onto $V$;
\item[(ii)] $H$ is polystable with pairwise non-isomorphic indecomposable direct summands.
\end{itemize}
\end{proposition}

\begin{proof}
By Lemma \ref{lemma2.1}, Proposition \ref{propa2} and Theorem \ref{thm3.3}, 
$E$ is generated by its sections, since $d \geq k > n$. By a standard result there is a surjective homomorphism 
$\cO^{n+1} \ra E$. This implies the existence of (\ref{eqn3.3}) and the assertion (i).

For the proof of (ii), note that every indecomposable direct summand of $H$ has negative degree. By dualizing the sequence (\ref{eqn3.3}),
we obtain a coherent system $(H^*,V^*)$ with $\dim V^* = k$. Moreover $(H^*,V^*)$ determines $(E,V)$. Hence (ii) follows by the same 
dimension counting argument as in the proof of Theorem \ref{thm3.3}.
\end{proof}

\section{Bounds for $\alpha$} 

\noindent
To obtain the range of $\alpha$ for which $G(\alpha;n,d,k) \neq \emptyset$, 
it is sufficient by Theorem \ref{cor3.4} to find the 
upper and lower bound for $\alpha$. For non-emptiness we must have 
$k \leq d$ or $(n,d,k)=(1,0,1)$, and by Proposition \ref{prop3.4} we can exclude 
the case $k=d$, gcd$(n,d) \neq 1$. We shall show that, in all other cases, $G(\alpha;n,d,k) \neq \emptyset$ for the full range of 
$\alpha$ permitted by \cite{bgn}. Since this has already been proved for $n=1$
in section 3, we shall suppose throughout that $n\ge2$. 

\begin{theorem} \label{thm4.1}
Suppose $k>0$ and either $k<d$ or $k=d$ and {\em gcd}$(n,d) =1$. Then 
\begin{center}
{\em inf} $\{\alpha \,|\,  G(\alpha;n,d,k) \neq \emptyset \} = 0.$
\end{center}  
\end{theorem}

\begin{proof}
We know that $\alpha > 0$ is a necessary condition for $G(\alpha;n,d,k)$ to be non-empty. It suffices to show that for small positive 
$\alpha$ there exist $\alpha$-stable coherent systems of type $(n,d,k)$ on $X$. Let  
$$
E = E_1 \oplus \cdots \oplus E_h
$$
be of rank $n$ and degree $d$ with $E_i$ pairwise non-isomorphic and stable of the same slope.

If $h=1$, then $E$ is stable, hence $(E,V)$ is $\alpha$-stable for small positive $\alpha$ for any choice of the subspace $V$.
We can therefore suppose $h>1$.
The only subbundles of $E$ which for small 
positive $\alpha$ might violate the $\alpha$-stability of $(E,V)$ are the subbundles of slope $\mu(E)$. But there are only 
finitely many of those, namely direct sums of some of the $E_j$. Hence, if $G$ is such a subbundle, it suffices to show 
that
\begin{eqnarray}
\frac{\dim (H^0(G) \cap V)}{\mbox{rk} (G)} < \frac{k}{n}  \label{eqn4.1}
\end{eqnarray}
for a general subspace $V$ of dimension $k$ of $H^0(E)$. By choosing $V$ generically we have, for all subbundles $G$ of slope
$\mu(E)$, 
$$
\dim (H^0(G) \cap V) = \max \{0,h^0(G) - \mbox{codim}(V)\}.
$$
If $\dim (H^0(G) \cap V) =0$, then (\ref{eqn4.1}) is certainly valid. So assume $h^0(G) - \mbox{codim}(V) \geq 0$, i.e. 
$\deg (G) = h^0(G) \geq$ codim$(V)$.
Then (\ref{eqn4.1}) becomes
$$
\frac{\deg(G) - \mbox{codim}(V)}{\mbox{rk}(G)} < \frac{d - \mbox{codim}(V)}{n},
$$
which is equivalent to $\mbox{rk}(G) < n$. This completes the proof of the theorem.
\end{proof}

We turn now to the question of the upper bound for the parameter $\alpha$, where we must distinguish the cases $k<n$ and $k \geq n$.

\begin{theorem} \label{thm4.2}
Suppose $0<k<n$. Suppose further that either $k<d$ or $k=d$ and $\gcd(n,d)=1$. Then 
\begin{center}
{\em sup} $\{ \alpha \,|\, G(\alpha;n,d,k) \neq 0\} = \frac{d}{n-k}$.
\end{center}
\end{theorem}
\begin{rem}
{\em Note that the statement of Theorem \ref{thm4.2} is covered by 
\cite[Remark 5.5]{bgn} if $\gcd(n-k,d) = 1$.}
\end{rem}

\begin{proof}[Proof of Theorem 5.2]
According to \cite[Lemma 4.1]{bgn} the number $\frac{d}{n-k}$ is an upper bound for $\alpha$. Hence it suffices to show 
that, for $\alpha$ slightly smaller than $\frac{d}{n-k}$, there exists an $\alpha$-stable coherent system of type $(n,d,k)$.

Let $(E,V)$ be as in Proposition \ref{prop3.5} and let $(E',V')$ be a proper coherent subsystem. We require to prove that 
\begin{eqnarray}
\mu_{\alpha}(E',V') < \mu_{\alpha}(E,V). \label{eqn4.2}
\end{eqnarray} 
The subbundle $E'$ determines an exact sequence
$$
0 \ra W' \ra E' \ra G' \ra 0,
$$
where $W'$ is a subbundle of $\cO^k$ and $G'$ is a subsheaf of the bundle $G$ of (\ref{eqn3.2}). The proof of 
\cite[Proposition 4.2]{bg} shows that (\ref{eqn4.2}) holds if $G'=0$ or $G'=G$ and also if $\mu(G') < \mu(G)$ (note 
that in the proof of \cite[Proposition 4.2]{bg} this is necessary in order to get $\delta > 0$). If $\mu(G') = \mu(G)$, 
the proof still works unless $W' \simeq \cO^{k'}$.

It remains to consider the following case:
$$
G = G_1 \oplus \cdots \oplus G_l
$$
with $G_i$ stable and pairwise non-isomorphic of the same slope, $G'$ a proper subbundle of $G$ with $\mu(G') = \mu(G)$
and $W' \simeq \cO^{k'}$ for some $k' \leq k$. Note that, in the case 
$k=d$ and $\gcd(n,d)=1$, $G$ is 
necessarily stable; so $G'$ does not exist. 

We may therefore assume that $k<d$.
Note that there are only finitely many such subbundles $G'$ of $G$. 
We now have an exact sequence 
$$
0 \ra \cO^{k'} \ra E' \ra G' \ra 0
$$
and we can suppose that $V' = H^0(\cO^{k'})$, i.e. dim$(V') = k'$. Note that $\mu(G') = \mu(G)$, so 
deg$(G') = \mbox{rk}(G')\cdot \frac{d}{n-k}$. Then we have
$$
\mu_{\alpha}(E',V')  - \mu_{\alpha}(E,V)  
 = \frac{\deg (G')}{k' + \mbox{rk}(G')} + \alpha \frac{k'}{k' + \mbox{rk}(G')} - \frac{d}{n} - \alpha \frac{k}{n} 
 $$
\begin{equation}
\qquad = \frac{k \cdot \mbox{rk}(G') -k'n +k'k}{n(n-k)(k'+\mbox{rk}(G'))}(d-\alpha(n-k)). \label{eqn4.3}
\end{equation}
We have to show that when $\alpha < \frac{d}{n-k}$, then, for a general extension 
\begin{equation}\label{ext}
0 \ra \cO^k \ra E \ra G \ra 0, 
\end{equation}
the right hand side of (\ref{eqn4.3}) is negative for all choices of $G'$ and $k'$. 
Equivalently we have to show 
$k \cdot \mbox{rk}(G') - k'n + k'k < 0$, i.e.
\begin{eqnarray}
k' > \frac{k \cdot \mbox{rk}(G')}{n-k}. \label{eqn4.4}
\end{eqnarray}

For the proof of (\ref{eqn4.4}) consider the diagram
$$
\xymatrix{
0 \ar[r] & \cO^k \ar[r] \ar@{=}[d] & E \ar[r] & G \ar[r] & 0\\
0 \ar[r] & \cO^k \ar[r] & E'' \ar[r] \ar[u] & G' \ar[r] \ar@{^{(}->}[u]_{i} & 0\\
0 \ar[r] & \cO^{k'} \ar[r] \ar@{^{(}->}[u] & E' \ar[r] \ar[u] & G' \ar[r] \ar@{=}[u] & 0.
  }              
$$
Let $(e_1,\ldots,e_k) \in H^1(G^*)^k$ denote the element classifying 
the extension (\ref{ext}). 
Then $(i^*e_1, \ldots, i^*e_k)$ classifies the pullback of (\ref{ext}) by the canonical embedding $i$.
The existence of the bottom row of the diagram implies that at most $k'$ of the elements $i^*e_1, \ldots, i^*e_k$ 
are linearly independent. Now, for the general extension (\ref{ext}), the number of these elements which are still linearly
independent is precisely $\min \{h^1(G'^*),k\}$. So 
$$
k' \geq \min \{h^1(G'^*),k\}.
$$
Since we have only finitely many subbundles $G'$, the general extension (\ref{ext}) has this property for all $G'$. Therefore we need 
to check (\ref{eqn4.4}) only for $k' = h^1(G'^*) = \deg (G')$ and $k' = k$. For $k'=k$ it is obvious and for $k' = \deg (G')$ it comes 
from the fact that $k<d$.
\end{proof}

\begin{theorem} \label{thm4.3}
Suppose $k \geq n$. Suppose further that either $k<d$ or $k=d$ and $\gcd(n,d)=1$. Then $G(\alpha;n,d,k) \neq \emptyset$ 
for arbitrarily large $\alpha$.
\end{theorem}

\begin{proof}
When $k=n$, this has already been proved in \cite[Theorem 5.6]{bgn}. So suppose $k>n$. For this we apply the results 
of \cite[section 6]{bgn}. Taking account of Theorem \ref{thm4.1}, it is enough to show that the dimension of the space of extensions
$$
0 \ra (E_1,V_1) \ra (E,V) \ra (E_2,V_2) \ra 0 
$$
is less than $\beta(d,k)$ whenever the following conditions hold:

\begin{itemize}
\item[(a)] $(E,V)$ is $\alpha$-stable of type $(n,d,k)$ for small positive $\alpha$;
\item[(b)] $(E_i,V_i)$ is of type $(n_i,d_i,k_i)$ for $i=1, 2$;
\item[(c)] $ \frac{k_1}{n_1} > \frac{k_2}{n_2}$ and $\frac{d_1}{n_1} < \frac{d_2}{n_2}$;
\item[(d)] $(E_i,V_i)$ is $\alpha$-stable for $\alpha$ slightly smaller than $\frac{n_1d_2-n_2d_1}{n_2k_1-n_1k_2}$ for $i=1, 2$.
\end{itemize}
Note that, since $k>n$, we have, by conditions (c) and (d), Proposition 
\ref{propa2} and Corollary \ref{cor3.2},
\begin{equation}\label{dkn}
n_1<k_1\le d_1.
\end{equation}

As in \cite[section 3]{bgn} we define 
$$
\HH^2_{21} = \mbox{Ext}^2((E_2,V_2),(E_1,V_1))
$$
and 
\begin{equation}
C_{12} = -d_2n_1 + d_1n_2+k_1d_2-k_1k_2. \label{eqn4.5}
\end{equation}
By \cite[Proposition 3.2]{bgn}, $\HH^2_{21} = H^0(E_1^* \otimes N_2)^*$, where $N_2$ is the kernel of the natural map 
$V_2 \otimes \cO \ra E_2$. 

By Lemma \ref{lemma3.1} every indecomposable summand of $E_1^*$ is of negative degree, 
while every indecomposable summand of $N_2$ has degree $\leq 0$. 
Since indecomposable bundles are semistable and we are in characteristic $0$,
it follows that $E_1^*\otimes N_2$ is a direct sum of semistable bundles
of negative degree; so $\HH^2_{21} = 0$.
Since also all moduli spaces of $\alpha$-stable coherent systems on elliptic curves have the correct dimensions by
Propositions \ref{propa1} and \ref{propa2} and Theorem \ref{thm3.3}, 
it follows from \cite[equation (18)]{bgn} that it is sufficient to prove that $C_{12} > 0$
under the conditions stated above. But
\begin{eqnarray*}
C_{12} & = & d_2(k_1-n_1)+d_1(n_2-k_2)+(d_1-k_1)k_2\\
& > & \frac{n_2d_1}{n_1}(k_1-n_1)+d_1(n_2-k_2)+(d_1-k_1)k_2\\
&& \left(\mbox{since} \;\; d_2 > \frac{n_2d_1}{n_1} \;\;  
\mbox{by assumption and} \; k_1>n_1 \; \mbox{by} \; (\ref{dkn})\right) \\
&=& d_1 \left( \frac{k_1n_2}{n_1} - k_2 \right) + (d_1-k_1)k_2\\
&> & 0\\
&& \left(\mbox{since} \; \frac{k_1}{n_1} > \frac{k_2}{n_2} \; \mbox{by assumption and} \, d_1 \geq k_1 \; 
\mbox{by}\; (\ref{dkn})\right).
\end{eqnarray*}
\end{proof}

\section{Variation of moduli spaces with $\alpha$}

\noindent
In this section, we again suppose throughout that $n\ge2$ and $k > 0$.
In the proof of Theorem \ref{thm4.3} we have used implicitly the idea of ``flips'' as defined in \cite{bgn}. These 
flips take place at the critical values $\alpha_i$ of the parameter $\alpha$, which have the form
$$
\alpha_i = \frac{n_1d_2-n_2d_1}{n_2k_1-n_1k_2},
$$
where $n_1+n_2=n,\; d_1+d_2=d$ and $k_1+k_2=k$. Moreover we can assume that
\begin{equation} 
\frac{k_1}{n_1} > \frac{k_2}{n_2} \quad \mbox{and} \quad \frac{d_1}{n_1} < \frac{d_2}{n_2}, \label{eqn5.1} 
\end{equation}
which together imply that $\alpha_i > 0$. If $k<n$, then we also require 
$$
\alpha_i < \frac{d}{n-k}
$$
(see Theorem \ref{thm4.2} or \cite[Proposition 4.2]{bgn}).
We can write this condition as 
\begin{equation}\label{alphai}
(n-k)\alpha_i < d,
\end{equation}
which makes it valid for all $k$.
For any $(n,d,k)$, there exist only finitely many critical values 
$$
0 < \alpha_1 < \alpha_2 < \cdots < \alpha_L
$$
(see \cite[Propositions 4.2 and 4.6]{bgn}).
For $\alpha_L< \alpha$ (or $\alpha_L<\alpha<\frac{d}{n-k}$ if $k<n$), 
the moduli space $G(\alpha;n,d,k)$ is independent of $\alpha$ and we denote it by $G_L$. 
Similarly, for $0 <\alpha < \alpha_1$ we denote $G(\alpha;n,d,k)$ by $G_0$.

\begin{lemma}\label{lemin}
Suppose that (\ref{eqn5.1}) and (\ref{alphai}) hold and 
$G(\alpha;n_1,d_1,k_1)$ and $G(\alpha;n_2,d_2,k_2)$ are both non-empty for
some $\alpha$. Then 
\begin{equation}
k_1 \leq d_1 \quad \mbox{and} \quad k_2 < d_2.  \label{eqn2}
\end{equation}
\end{lemma}


\begin{proof} By (\ref{eqn5.1}), $k_1\ge1$; so, by Proposition \ref{propa2} and Corollary \ref{cor3.2}, $k_1>d_1$
only when $(n_1,d_1,k_1)=(1,0,1)$. In this case, $k_2<n_2$ by (\ref{eqn5.1}), so $k<n$. Moreover $\alpha_i=\frac{d}{n-k}$,
which contradicts (\ref{alphai}). Hence $k_1\le d_1$. Together with 
(\ref{eqn5.1}) this implies $\frac{k_2}{n_2} < \frac{k_1}{n_1} \leq \frac{d_1}{n_1} < \frac{d_2}{n_2}$, which
completes the proof.
\end{proof}

As $\alpha$ increases through $\alpha_i$, we must delete from the moduli space $G(\alpha;n,d,k)$ a closed subset 
$G_i^-$ and insert $G_i^+$. In \cite[section 6]{bgn} an explicit description of $G_i^-$ and $G_i^+$ is given together
with estimates for their codimensions. In our case, since all moduli spaces
have the expected dimensions and $\HH^2_{21}=\HH^2_{12}=0$, we have
$$
\mbox{codim} (G_i^-) \ge \min C_{12} \quad \mbox{and} 
\quad \mbox{codim} (G_i^+) \ge \min C_{21},
$$
where $C_{12}$ is given by (\ref{eqn4.5}) and 
$$
C_{21} = -d_1n_2+d_2n_1+k_2d_1-k_1k_2
$$
and the minima are taken over all values of $n_1$, $d_1$, $k_1$ satisfying $(\ref{eqn5.1})$ and giving rise to the critical 
value $\alpha_i$.

\begin{proposition} \label{prop1}
Suppose that (\ref{eqn5.1}) and (\ref{alphai}) hold and 
$G(\alpha;n_1,d_1,k_1)$ and $G(\alpha;n_2,d_2,k_2)$ are both non-empty for
some $\alpha$. Then both $C_{12}$ and $C_{21}$ are positive.
\end{proposition}

\begin{proof}
$C_{21} = n_1n_2(\frac{d_2}{n_2} - \frac{d_1}{n_1}) + k_2(d_1 - k_1) >0$ by (\ref{eqn5.1}) and (\ref{eqn2}).

As for $C_{12}>0$, this is proved in the course of the proof of Theorem \ref{thm4.3} for $k>n$ and, in fact, the proof is valid for 
$k \geq n$. For $k<n$ we have, by (\ref{alphai}), 
$$
\frac{n_1d_2-n_2d_1}{n_2k_1-n_1k_2} < \frac{d}{n-k}
$$ and thus
$$
\begin{array}{ll}
n(n_1d_2-n_2d_1) & < k(n_1d_2-n_2d_1)+d(n_2k_1-n_1k_2)\\
& = n(k_1d_2-k_2d_1).
\end{array} 
$$
Therefore
$$
\begin{array}{ll}
0 & < d_1n_2-d_2n_1+k_1d_2-k_1k_2+k_1k_2-d_1k_2\\
& = C_{12} + k_2(k_1-d_1)
\end{array}
$$
and hence $C_{12} > k_2(d_1-k_1) \geq 0$.
\end{proof}

\begin{rem} \label{rem2}
{\em Proposition \ref{prop1} shows that all flips are ``good'' in the sense of \cite{bgn}.
It also gives an alternative proof of Theorem \ref{thm4.2}, but the proof in section 5 gives more information 
about the structure of $G_L$.}
\end{rem}

\begin{ex} \label{ex5.3}
{\em Suppose $k=d$. For non-emptiness of $G(\alpha;n,d,k)$ we require gcd$(n,d) = 1$. The inequalities 
(\ref{eqn2}) have no common solution. Hence there are no critical values and $G_0=G_L = X$. 
In fact, $(E,V)$ is $\alpha$-stable for some $\alpha$ if and only if $E$ is a stable bundle (and $V=H^0(E))$.}
\end{ex}

Suppose now $k=d-1$. According to Lemma \ref{lemin} we must have 
$$
k_1 = d_1 \quad \mbox{and} \quad k_2=d_2-1.
$$
(\ref{eqn5.1}) then gives 
\begin{equation}
\frac{d_2}{n_2} > \frac{d_1}{n_1} > \frac{d_2-1}{n_2}, \label{eqn3}
\end{equation}
i.e. 
$$
n_2d < nd_2<n_2d+n-n_2.
$$
When $d \leq n$, we have also $\alpha_i < \frac{d}{n-k}$, which is equivalent to
$$
C_{12} > k_2(d_1-k_1) = 0. 
$$
This in turn is equivalent to
\begin{equation}
d_1n_2-d_2n_1+d_1> 0.   \label{eqn4}
\end{equation}
Since in this case $d_1<n_1$, this is strictly stronger than the second inequality of (\ref{eqn3}). 
Note that $d_1 = 0$ is not allowed because of (\ref{eqn3}).

\begin{ex} \label{ex5.4}
{\em Let $d=2, k=1$. According to our discussion above the only possibility is $k_1=d_1=1,\; d_2=1,\; k_2=0$.
But this contradicts (\ref{eqn3}) and (\ref{eqn4}). So there are no critical values and $G_0=G_L$.

If $n$ is odd, then $(E,V)$ is $\alpha$-stable for small $\alpha$ if and only if $E$ is stable as a bundle and $V$ 
is any one-dimensional subspace of the two-dimensional space $H^0(E)$. So $G_0$ is a $\PP^1$-bundle over $X$.

If $n$ is even, we can use the standard description of $G_L$ (see \cite[Remark 5.5]{bgn}). The points of $G_L$ correspond bijectively to 
the classes of non-trivial exact sequences 
$$
0 \ra \cO \ra E \ra G \ra 0
$$
with $G$ stable of rank $n-1$ and degree 2. So again $G_0=G_L$ 
is a $\PP^1$-bundle over $X$.}
\end{ex} 

\begin{ex} \label{ex5.5}
{\em Let $d=3,\; k=2$. The only possible solutions of (\ref{eqn3}) and (\ref{eqn4}) are $d_1=2$ and $d_2=1$, implying 
$k_1=2, \;k_2=0$. (\ref{eqn3}) and (\ref{eqn4}) imply
$$n_1> 2n_2> n_1-2,
$$
i.e. $n_1=2n_2+1$. Thus 
$$
n = n_1+n_2 = 3n_2+1 \equiv 1\; \mbox{mod} \;3.
$$
Hence, if $n\not\equiv 1 \; \mbox{mod} \; 3$, there are no critical values. Then the same argument as in 
the previous example shows that $G_0=G_L$ is a $\PP^2$-bundle over $X$.

Finally, assume $n \equiv 1 \; \mbox{mod} \; 3$. Then both gcd$(n,d)=1$ and gcd$(n-k,d) =1$, so $G_0$ and $G_L$ are both 
$\PP^2$-bundles over $X$. However there is one critical value $\alpha_1 = \frac{1}{2n_2}$. Moreover $C_{12}=C_{21} = 1$. 
So $G_1^-$ and $G_1^+$ both have dimension at most 2.

Now $G_1^-$ is given by non-trivial extensions
\begin{equation} 
0 \ra (E_1,V_1) \ra (E,V) \ra (E_2,V_2) \ra 0,  \label{eqn5}
\end{equation}
where $(E_1,V_1)$ is $\alpha_1$-stable of type $(2n_2+1,2,2)$, $E_2$ is a stable bundle of type $(n_2,1) $ and $V_2 =0$.
According to \cite[equation (8)]{bgn} and the fact that $\HH^0_{21} =\HH^2_{21} = 0$, we have that
$$
\dim(\mbox{Ext}^1((E_2,V_2),(E_1,V_1))) = C_{21}=1.
$$
So there is a unique non-trivial extension for any $(E_1,V_1)$ and $(E_2,V_2)$, implying that
$$
G_1^- \simeq X \times X.
$$ 
Similarly, $G_1^+$ is given by non-trivial extensions 
\begin{equation} 
0 \ra (E_2,V_2) \ra (E,V) \ra (E_1,V_1) \ra 0  \label{eqn6}
\end{equation}
with $(E_1,V_1),\; (E_2,V_2)$ as above. Thus
$$
\dim (\mbox{Ext}^1((E_1,V_1),(E_2,V_2))) = C_{12}= 1.
$$
Again there is a unique non-trivial extension for any $(E_1,V_1)$ and $(E_2,V_2)$, implying that
$$
G_1^+ \simeq X \times X.
$$ 

As $\alpha \ra \alpha_1$ from below, $G_0$ ``degenerates'' to a moduli space $\tilde{G}(\alpha_1)$ of S-equivalence 
classes of $\alpha_1$-semistable coherent systems. In fact, the coherent system $(E,V) \in G_0$ given by (\ref{eqn5})
is replaced by the S-equivalence class defined by $(E_1,V_1) \oplus (E_2,V_2)$.
Similarly, as $\alpha \ra \alpha_1$ from above, $G_L$ ``degenerates'' to $\tilde{G}(\alpha_1)$ and the coherent system 
$(E,V) \in G_L$ given by (\ref{eqn6}) is replaced by the same S-equivalence class.
It is possible that $\tilde{G}(\alpha_1)$ has singularities and may even be
reducible, but the natural morphisms $G_0\to\tilde{G}(\alpha_1)$ and
$G_L\to\tilde{G}(\alpha_1)$ are both bijective onto the same irreducible
component of $\tilde{G}(\alpha_1)$. Since $G_0$ and $G_L$ are smooth, they are
both normalisations of this component and are therefore isomorphic.
However the families parametrised by $G_0$ and $G_L$ are different:
the coherent systems given by (\ref{eqn5}) belong to $G_0$ and 
are replaced by the coherent systems given by 
(\ref{eqn6}), which belong to $G_L$.}
\end{ex}

\section{Poincar\'e and Picard bundles}

\noindent
We need some facts about Poincar\'e and Picard bundles for elliptic curves. 
For $n=2$, the relevant facts are contained in \cite[Calculation 5.14]{bdw},
but for higher rank we cannot find a reference, 
although the situation for general genus is described in \cite{ab}. The main difference is that in our case the 
moduli spaces are identified with the curve $X$ and we need to take this into account.

Fix $n$ and $d>0$ with gcd$(n,d) = 1$ and let $M = M(n,d)$ denote the moduli space of stable bundles on $X$ 
of rank $n$ and degree $d$. By \cite[Theorem 7]{at} the map
$$
\det: M \ra X
$$
is an isomorphism. 
Let $p_1$ and $p_2$ denote the projections of $X \times M$ and $q_1$ and $q_2$ the projections of $X \times X$.
Hence, if $\cP$ denotes the Poincar\'e bundle parametrizing line bundles of degree $d$ on $X \times X$ and
$\cU$ is a Poincar\'e bundle on $X \times M$, then 
$$
\det(\cU) \simeq p_2^*L_1 \otimes (\mbox{id} \times {\det})^*(\cP)
$$ 
with $L_1 \in \mbox{Pic}(M)$. By \cite[p. 335]{acgh}, 
$$
c_1(\cP) = d[X] + \xi_1
$$ with $[X] = [q_1^*(\mbox{point})]$ and $\xi_1^2= -2[X \times X]$.
So
$$
c_1(\cU) = d[X] + \xi_1 + p_2^*(a_1),
$$
where $a_1 = c_1(L_1)$, $[X] = [p_1^*(\mbox{point})]$ and, by abuse of notation, we write $\xi_1 = (\mbox{id} \times \det)^*\xi_1$. 
So now
$$
\xi_1^2 = -2[X \times M].
$$
We write
$$
c_2(\cU) = p_2^*(f_2)[X]
$$
with $f_2 \in H^2(M)$. Note that ${p_2}_*\cU$ is a vector bundle of rank $d$ on $M$, usually called the Picard bundle.

\begin{proposition} \label{prop5.1}
Let $[M]$ denote the fundamental class of $M$. Then
\begin{itemize}
\item[(i)] $ \;c_1({p_2}_*\cU) = da_1 - [M] - f_2$;
\item[(ii)] $ \; nf_2 = (n-1)(da_1 - [M])$.
\end{itemize}\end{proposition}

\begin{proof}
(i): Since $p_2^*(a_1)\xi_1 = 0$ and $[X]\xi_1 = 0$,
$$
ch_2(\cU) = \frac{1}{2}c_1^2(\cU) - c_2(\cU) = dp_2^*(a_1)[X] 
- [X \times M] - p_2^*(f_2)[X].
$$ 
Since $H^1(E) = 0$ for all $E \in M$, Grothendieck-Riemann-Roch gives
$c_1({p_2}_*\cU) = da_1 - [M] -f_2$ as required.\\

\noindent(ii): Here we apply Grothendieck-Riemann-Roch to the bundle $End(\cU)$.
$$
\begin{array}{cl}
c_2(End (\cU)) & = -(n-1)c_1^2(\cU) + 2nc_2(\cU)\\
& = (n-1)\{2[X \times M] - 2dp_2^*(a_1)[X]\} +2np_2^*(f_2)[X].
\end{array}
$$
Since $c_1(End(\cU)) =0$, we get
$$
ch_2((End(\cU)) = (n-1)\{2dp_2^*(a_1)[X]-2[X \times M]\} - 2np_2^*(f_2)[X].
$$
Hence
$$
c_1({p_2}_*End(\cU)) - c_1(R^1{p_2}_*(End(\cU))) = 2(n-1)(da_1-[M]) - 2nf_2.
$$
But ${p_2}_*End(\cU) \simeq R^1{p_2}_*(End(\cU))^*$ by relative Serre duality.
Moreover $R^1{p_2}_*(End(\cU))$ is the tangent bundle of $M$, which is trivial. This gives the result.
\end{proof}

Since $H^2(M) \simeq \ZZ$, we can write
$$
a_1 = r[M].
$$ Then Proposition \ref{prop5.1}(ii) reads $nf_2 = (n-1)(dr - 1)[M]$ and we deduce
$$
dr-1=sn \quad \mbox{and} \quad f_2=(n-1)s[M]
$$
for some integer $s$. 

\begin{proposition} \label{prop5.2}
$ \;\; c_1({p_2}_*\cU) = s[M]$.
\end{proposition}

\begin{proof}
$ c_1({p_2}_*\cU) = dr[M] - [M] - (n-1)s[M] = s[M]$.
\end{proof}

\begin{rem} \label{rem5.3}
{\em Write $d_1=c_1({p_2}_*\cU)$. Since $rd-sn =1$,
$$
[M] = da_1 - nd_1.
$$
So $a_1$ and $d_1$ generate $H^2(M)$. Although this is in line with \cite[Theorem 9.11]{ab}, it does not seem 
to follow directly from that theorem}.
\end{rem}

\section{An Example}

\noindent
Suppose $0<k<n, \; k<d$ with gcd$(n,d) = 1$ and gcd$(n-k,d) = 1$.

\begin{proposition} \label{prop7.7}
Under the above assumptions,
\begin{itemize}
\item[(i)] $G_0$ can be identified with the Grassmannian bundle over $M=M(n,d)$ of subspaces of dimension $k$ 
in the fibres of the Picard bundle $\cE = {p_2}_*\cU$;
\item[(ii)] $G_L$ can be identified with the Grassmannian bundle over $M(n-k,d)$ of subspaces of dimension $k$ in the fibres of the 
Picard bundle $\cF = ({p_2}_*\cU')^*$, where $\cU'$ is a Poincar\'e bundle on $X \times M(n-k,d)$;
\item[(iii)] $c_1(\cF) = - s'[M(n-k,d)]$ for some $s'$ such that $s'(n-k) \equiv -1$ {\em mod} $d$.  
\end{itemize}
\end{proposition}

\begin{proof}
(i): Since gcd$(n,d) = 1$, 
$$
(E,V) \in G_0 \quad \Leftrightarrow \quad E \quad \mbox{is stable}.
$$
The result follows immediately.

\noindent(ii): Since gcd$(n-k,d) = 1$, it follows from \cite[Remark 5.5]{bgn} that $G_L$ is a Grassmannian 
bundle associated with the vector bundle $R^1{p_2}_*(\cU'^*)$. By relative Serre duality this is isomorphic to $({p_2}_*\cU')^*$.

\noindent(iii) follows at once from Proposition \ref{prop5.2} applied to $\cU'$.
\end{proof}

After identifying $M$ and $M(n-k,d)$ with $X$ via the determinants, we can regard $\cE$ and $\cF$ as bundles over 
the elliptic curve $X$. 

\begin{proposition} \label{prop6.1}
The projective bundles $P(\cE)$ and $P(\cF)$ are independent of the choices of the Poincar\'e bundles $\cU$ and $\cU'$.
Moreover, if 
$$
s + s' \not\equiv 0 \; \mbox{mod} \; d,
$$
then $P(\cE)$ and $P(\cF)$ are not isomorphic as varieties.
\end{proposition}

\begin{proof}
The Poincar\'e bundles $\cU$ and $\cU'$ are determined up to tensoring by a line bundle on the moduli spaces. So
the same holds for $\cE$ and $\cF$ by the projection formula. Hence $P(\cE)$ and $P(\cF)$ are independent of the 
choices of $\cU$ and $\cU'$. Since $c_1(\cE) = s[X],\, c_1(\cF) = -s'[X]$ and $\cE$ and $\cF$ both have rank $d$, it follows
that, if $s + s' \not\equiv 0 \; \mbox{mod}\;  d$, then $P(\cE)$ and $P(\cF)$ are not isomorphic as projective bundles over $X$.
Now any isomorphism of varieties from $P(\cE)$ to $P(\cF)$ must map fibres to fibres and induces a linear map on each fibre.
Hence there exists an automorphism $\psi$ of $X$ such that $P(\cE)$ and $P(\psi^*\cF)$ are isomorphic as projective 
bundles over $X$. Since $c_1(\psi^*\cF) = c_1(\cF)$, this proves the last part of the proposition.
\end{proof}

From the proposition we can derive examples showing that the moduli spaces $G_0$ and $G_L$ can be non-isomorphic. 
We finish with a specific example.

\begin{ex}
{\em If $k = 1$, we can identify $G_0$ with $P(\cE)$ and $G_L$ with $P(\cF)$ by 
Proposition \ref{prop7.7}. Moreover $s$ and $s'$ are determined modulo $d$ by the conditions
$$
sn \equiv -1 \; \mbox{mod} \;d \quad \mbox{and} \quad s'(n-1) \equiv -1\; \mbox{mod} \; d.
$$
If  $s+s' \equiv 0 \; \mbox{mod} \; d$, adding these conditions gives
$$
s \equiv -2 \; \mbox{mod} \; d \quad \mbox{and} \quad s' \equiv +2 \; \mbox{mod} \; d.
$$
Hence, if we choose an example for the pair $(n,d)$ such that
$$
2n \not\equiv 1 \; \mbox{mod} \; d,
$$
then the moduli spaces $G_0 \simeq P(\cE)$ and $G_L \simeq P(\cF)$ for the triple $(n,d,1)$ are non-isomorphic. 
For an example, choose $n=5$ and 
$d=7$.}
\end{ex}

\end{document}